\documentclass{amsart}

\usepackage{amsmath,amsthm,amssymb}

\newtheorem{theorem}{Theorem}

\numberwithin{theorem}{section}

\theoremstyle{remark}
\newtheorem{remark}[theorem]{Remark}

\begin{document}

\title{Further Rigid Triples of Classes in {$G_{2}$}}
\author[M.\ Conder]{Matthew Conder}
\address{M.J.\ Conder, Department of Pure Mathematics and Mathematical Statistics, Centre for Mathematical Sciences, University of Cambridge, Wilberforce Road, Cambridge, CB3 0WB, United Kingdom}
\email{mjc271@cam.ac.uk}
\thanks{The first author is jointly supported by the Cambridge and Woolf Fisher Trusts.}

\author[A.\ Litterick]{Alastair Litterick}
\address{A.J.\ Litterick, Fakult\"{a}t f\"{u}r Mathematik, Ruhr-Universit\"{a}t Bochum, Universit\"{a}tsstra{\ss}e 150, D-44780 Bochum, Germany \and Fakult\"{a}t f\"{u}r Mathematik, Universit\"{a}t Bielefeld, Postfach 100131, D-33501 Bielefeld, Germany}
\email{ajlitterick@gmail.com}
\thanks{The second author is supported by the Alexander von Humboldt Foundation.}



\begin{abstract}
We establish the existence of two rigid triples of conjugacy classes in the algebraic group $G_{2}$ in characteristic $5$, complementing results of the second author with Liebeck and Marion. As a corollary, the finite groups $G_{2}(5^n)$ are not $(2,4,5)$-generated, confirming a conjecture of Marion in this case.
\end{abstract}

\maketitle

\section{Introduction}

Let $G$ be a connected simple algebraic group over an algebraically closed field $K$, and let $C_{1}$, $\ldots$, $C_{s}$ be conjugacy classes of $G$. Following \cite{Strambach1999}, we say the $s$-tuple $\mathbf{C} = (C_{1},\ldots,C_{s})$ is \emph{rigid} in $G$ if the set
\[ \mathbf{C}_{0} \stackrel{\textup{def}}{=} \{ (x_{1},\ldots,x_{s}) \in C_{1} \times \ldots \times C_{s} \, : \, x_{1}x_{2}\ldots x_{s} = 1 \} \]
is non-empty and forms a single orbit under the action of $G$ by simultaneous conjugation.

Some well-known examples of rigid tuples of classes in simple algebraic groups are the Belyi triples and Thompson tuples, defined in \cite{Voelklein1998}. Other rigid triples are known, see for instance \cite{Dettweiler1999,Feit1985,Guralnick2014,Liebeck2011,Thompson1985}. Rigid tuples of classes are interesting in the context of the inverse Galois problem \cite{Malle1999}, and also arise naturally in the theory of ordinary differential equations \cite{Katz1996}.

Recall that a group is $(a,b,c)$-generated if it is generated by elements $x$, $y$ and $z$, of respective orders $a$, $b$ and $c$, such that $xyz = 1$. The group is then called an $(a,b,c)$-group, and the triple $(x,y,z)$ is called an $(a,b,c)$-triple of the group. The theory of $(a,b,c)$-generation of finite groups has close connections to rigidity, for instance it is a basic observation that given a rigid tuple $\mathbf{C}$ of classes of $G$, all subgroups $\left<x_1,\ldots,x_s\right>$ for $(x_1,\ldots,x_s) \in \mathbf{C}_{0}$ are conjugate in $G$, so that there is at most one $r > 0$ such that the finite subgroup $G(p^{r})$ is generated by elements in such an $s$-tuple.

Let $K = \bar{\mathbb{F}}_{5}$ be the algebraic closure of the field of five elements. In \cite{Liebeck2011} it is shown that the simple algebraic group $G = G_{2}(K)$ has a rigid triple of conjugacy classes of elements of orders $2$, $5$ and $5$, and any triple of elements $(x_1,x_2,x_3)$ in the corresponding set $\mathbf{C}_{0}$ generates a copy of $\textup{Alt}_{5}$. This is then used to show that none of the groups $G_{2}(5^n)$, $SL_{3}(5^n)$ or $SU_{3}(5^n)$ is a $(2,5,5)$-group.

Here we produce two further rigid triples of classes in $G = G_2(K)$, closely related to the triple above. Recall from \cite{Chang1974} that $G$ has a unique class of involutions, with representative $t$, say, and $C_{G}(t) = A_{1}\tilde{A}_{1}$ is a central product of two subgroups $SL_{2}(K)$, where $A_{1}$ (resp.\ $\tilde{A}_{1}$) is generated by a long (resp.\ short) root subgroup of $G$. There also exist two classes of elements of order $4$, with representatives $s_{1}$ and $s_{2}$, such that $C_{G}(s_{1}) = A_{1}T'$ and $C_{G}(s_{2}) = \tilde{A}_{1}T''$, where $T'$ and $T''$ are $1$-dimensional tori. Finally, recall from \cite{Lawther1995} that $G$ has three classes of unipotent elements of order 5: the long and short root elements, and the class labelled $G_{2}(a_1)$, with representative $u = x_{\beta}(1)x_{3\alpha + \beta}(1)$, where $\alpha$ (resp.\ $\beta$) is the short (resp.\ long) simple root of $G$. From  \cite[Table 22.1.5]{Liebeck2012}, the centraliser $C_{G}(u) = U_{4}.\textup{Sym}_{3}$, where $U_{4}$ is a $4$-dimensional connected unipotent group.

\begin{theorem} \label{thm:main} \leavevmode
\begin{itemize}
\item[\textup{(i)}]The triples of classes $\mathbf{C} = (t^{G},s_{1}^{G},u^{G})$ and $\mathbf{D} = (t^{G},s_{2}^{G},u^{G})$ are rigid in $G = G_{2}(K)$.
\item[\textup{(ii)}]Every triple of elements in $\mathbf{C}_{0}$ or $\mathbf{D}_{0}$ generates a subgroup isomorphic to the symmetric group $\textup{Sym}_{5}$.
\item[\textup{(iii)}]None of the groups $G_{2}(5^n)$ are a $(2,4,5)$-group for any $n$. Neither are the groups $SL_{3}(5^n)$ or $SU_{3}(5^n)$.
\end{itemize}
\end{theorem}

\begin{remark} \leavevmode
\begin{enumerate}
\item Each subgroup $\textup{Sym}_{5}$ in part (ii) here contains a subgroup $\textup{Alt}_{5}$ arising from \cite[Theorem 1(ii)]{Liebeck2011}.
\item Keeping track of details in the proof in \cite{Liebeck2011} shows that $G_{2}(K)$ has a unique class of subgroups $\textup{Alt}_{5}$. These subgroups have centraliser $\textup{Sym}_{3}$, and by Lang's theorem these split into three classes in $G_{2}(5^r)$, with centraliser orders $6$, $3$ and $2$. Similarly, if $S$ and $S'$ are representatives of the two subgroup classes in part (ii) here, then $C_{G}(S) \cong \textup{Sym}_{3}$, while $C_{G}(S')$ is cyclic of order $2$. It follows that the class of $S$ (resp.\ $S'$) splits into $3$ (resp.\ $2$) classes of subgroups in $G_{2}(5^r)$, with centralisers of order $6$, $3$, $2$ (resp.\ $2$ and $2$).
\item A conjecture of Marion \cite{Marion2010} states that, for a simple algebraic group $G$ in characteristic $p$, if $\delta_{i}$ denotes the dimension of the variety of elements of $G$ of order $i$ and if $\delta_{a} + \delta_{b} + \delta_{c} = 2\, \textup{dim}(G)$, then at most finitely many of the finite groups $G(p^{r})$ are $(a,b,c)$-groups. For $G$ of type $G_{2}$, this criterion holds precisely when $(a,b,c) = (2,4,5)$ or $(2,5,5)$. Hence part (iii), together with \cite[Theorem 1(iii)]{Liebeck2011}, verifies the conjecture for $G = G_{2}$ in characteristic $5$. A non-constructive proof of this fact is given in \cite[Proposition 3.7(i)]{Jambor2017}, where it is shown that every $(2,4,5)$-subgroup and $(2,5,5)$-subgroup of $G_{2}(K)$ is reducible on the natural $7$-dimensional module, by considering the dimensions of $SL_{7}(K)$-conjugacy classes of elements in the relevant $(a,b,c)$-triples.
\end{enumerate}
\end{remark}

\section{Proof of the Theorem}

We proceed in the manner of \cite{Liebeck2011}. Let $G = G_{2}(K)$ and $t$, $u$, $s_{1}$, $s_{2} \in G$ as above. If $\sigma$ is a Frobenius morphism of $G$ induced from the field map $x \mapsto x^{5}$ of $K$, then
\[ G = \bigcup_{n=1}^{\infty} G_{\sigma^{n}} = \bigcup_{n = 1}^{\infty} G_{2}(5^n). \]
The element $u = x_{\beta}(1)x_{3\alpha + \beta}(1)$ is a regular unipotent element in a subgroup $A_{2} = SL_{3}(K)$ of $G$ generated by long root groups, and therefore lies in a subgroup $\Omega_{3}(5) \cong \textup{Alt}_{5}$ of $G$, which we denote by $A$. Now, let $S = N_{A_2}(A) = SO_{3}(5) \cong \textup{Sym}_{5}$. Following the proof given in \cite{Liebeck2011} we find that $N_{G}(A) = S \times C_{G}(A)$ and $C_{G}(A) = \left<z,\tau\right> \cong \textup{Sym}_{3}$, where $\left<z\right>$ is the centre of $A_{2}$ and $\tau$ is an outer involution in $N_{G}(A_{2}) = A_{2}.2$. Note that $C_{A_2}(\tau) = SO_{3}(K)$, so $\tau \in C_{G}(S)$.

Let $v$ be an involution in $S \setminus A$, so that $S = \left<A,v\right>$, and define $S' = \left<A,v\tau\right>$, so that $S' \cong \textup{Sym}_{5}$ also. Then $C_{G}(S)$, $C_{G}(S') \le C_{G}(A) = \left<z,\tau\right>$ and therefore
\begin{align}
C_{G}(S) &= \left<z,\tau\right>, \label{cgs} \\
C_{G}(S') &= \left<\tau\right>. \label{cgs'}
\end{align}
In particular $S$ and $S'$ are not conjugate in $G$.

Next consider the set of $(2,4,5)$-triples of $\textup{Sym}_{5}$. It is straightforward to show that there are exactly $120$ such triples, and that $\textup{Sym}_{5}$ acts transitively on these by simultaneous conjugation.

Now let $\mathbf{C} = (t^{G},s_1^{G},u^{G})$ and $\mathbf{D} = (t^{G},s_{2}^{G},u^{G})$, and for $q$ a fixed power of $5$ let $\mathbf{C}_{0}(q) = \mathbf{C}_{0} \cap G_{2}(q)^{3}$ and $\mathbf{D}_{0}(q) = \mathbf{D}_{0} \cap G_{2}(q)^{3}$. We now show that $|\mathbf{C}_{0}(q)| = |\mathbf{D}_{0}(q)| = |G_{2}(q)|$. For this we require the character table of $G_{2}(q)$, given in \cite{Chang1974} and available in the CHEVIE \cite{Geck1996} computational package. Since $C_{G}(u)/C_{G}(u)^{\circ} = S_{3}$, an application of Lang's theorem \cite[Theorem 21.11]{Malle2011} shows that $u^{G} \cap G_{2}(q)$ splits into three classes, with representatives denoted in \cite{Chang1974} by $u_{3}$, $u_{4}$ and $u_{5}$, and respective centraliser orders $6q^{4}$, $3q^{4}$ and $2q^{4}$. For $x,y,z \in G_{2}(q)$ let $a_{xyz}$ be the corresponding class algebra constant. Calculations with the character table show that
\begin{align*}
a_{ts_{i}u_{j}} &= \left\{
\begin{array}{cl}
q^{4} & \textup{ if } i = 1,\ j \in \{3,4,5\} \textup{ or } i = 2,\ j = 4, \\
3q^{4} & \textup{ if } i = 2,\ j = 3, \\
0 & \textup{ if } i = 2, j = 5.
\end{array}
\right.
\end{align*}
and it follows that 
\begin{align*}
|\mathbf{C}_{0}(q)| &= \sum_{j = 3}^{5} |u_{j}^{G_2(q)}|a_{ts_{1}u_{j}} = |G_{2}(q)| \left(\frac{q^{4}}{6q^{4}} + \frac{q^{4}}{3q^{4}} + \frac{q^{4}}{2q^{4}}\right) = |G_{2}(q)|, \\
|\mathbf{D}_{0}(q)| &= \sum_{j = 3}^{5} |u_{j}^{G_2(q)}|a_{ts_{2}u_{j}} = |G_{2}(q)| \left(\frac{3q^{4}}{6q^{4}} + \frac{q^{4}}{2q^{4}}\right) = |G_{2}(q)|.
\end{align*}

Now let $\mathbf{E}$ denote (resp.\ $\mathbf{E}'$) denote the set of triples $(x_1,x_2,x_3) \in \mathbf{C}_{0} \cup \mathbf{D}_{0}$ which generate a conjugate of $S$ (resp.\ a conjugate of $S'$). Then $G$ is transitive on both $\mathbf{E}$ and $\mathbf{E}'$, since if $\left<x_1,x_2,x_3\right> = \left<y_1,y_2,y_3\right>^{g}$ are each isomorphic to $\textup{Sym}_{5}$, then $(x_1^{g},x_2^{g},x_3^{g})$ and $(y_1,y_2,y_3)$ are $(2,4,5)$ triples in a fixed copy of $\textup{Sym}_{5}$, hence conjugate in $\textup{Sym}_{5}$ by the observation above. Moreover both $\mathbf{E}$ and $\mathbf{E}'$ are non-empty, since $S$ and $S'$ each contain $(2,4,5)$-triples and a unique conjugacy class of unipotent elements, whose elements are conjugate to an element of $A$ and therefore are conjugate to $u$. By (\ref{cgs}) and (\ref{cgs'}) the stabiliser of a point in $\mathbf{E}$ is isomorphic to $\textup{Sym}_{3}$, and the stabiliser of a point in $\mathbf{E}'$ is cyclic of order $2$. Hence applying Lang's theorem shows that $\mathbf{E}(q) = \mathbf{E} \cap G_{2}(q)^{3}$ splits into three $G_{2}(q)$-orbits, of orders $|G_{2}(q)|/r$ for $r = 2,3,6$, and similarly $\mathbf{E}'(q) = \mathbf{E}' \cap G_{2}(q)^{3}$ splits into two orbits, each of order $|G_{2}(q)|/2$. Therefore,
\begin{align*}
|\mathbf{E}(q)| + |\mathbf{E}'(q)| = |G_{2}(q)|\left(\frac{1}{6} + \frac{1}{3} + \frac{1}{2} + \frac{1}{2} + \frac{1}{2}\right) = |\mathbf{C}_{0}(q)| + |\mathbf{D}_{0}(q)|
\end{align*}
and it follows that $\mathbf{C}_{0}(q) \cup \mathbf{D}_{0}(q) = \mathbf{E}(q) \cup \mathbf{E}'(q)$ for each $q$. Therefore
\[ \mathbf{C}_{0} \cup \mathbf{D}_{0} = \bigcup_{n = 1}^{\infty}\mathbf{C}_{0}(5^n) \cup \mathbf{D}_{0}(5^n) = \bigcup_{n = 1}^{\infty}\mathbf{E}(5^n) \cup \mathbf{E}'(5^n) = \mathbf{E} \cup \mathbf{E}' \]
Hence $G$ has exactly two orbits on $\mathbf{C}_{0} \cup \mathbf{D}_{0}$. A triple in $\mathbf{C}_{0}$ cannot lie in the same orbit as a triple in $\mathbf{D}_{0}$ since the corresponding elements of order $4$ are not $G$-conjugate, and it follows that the two $G$-orbits are $\mathbf{C}_{0}$ and $\mathbf{D}_{0}$.

This proves parts (i) and (ii) of the Theorem. For part (iii), suppose that $G_{2}(5^n)$, $SL_{3}(5^n)$ or $SU_{3}(5^n)$ is a $(2,4,5)$-group, with corresponding set of generators $x_1$, $x_2$, $x_3$. Since $L(G_2) \downarrow A_{2}$ is a direct sum of $L(A_2)$ and two 3-dimensional irreducible $A_{2}$-modules (cf.\ \cite[Table 8.5]{Liebeck1996a}), it follows that $C_{L(G_2)}(x_1,x_2,x_3) = 0$. An application of a result of Scott \cite{Scott1977} to the module $L(G)$, as in the proof of \cite[Corollary 3.2]{Strambach1999}, then yields
\[ \textup{dim}(x_{1}^{G}) + \textup{dim}(x_{2}^{G}) + \textup{dim}(x_{3}^{G}) \ge 2\, \textup{dim}(G) = 28, \]
implying $(x_{1}^{G},x_{2}^{G},x_{3}^{G}) = \mathbf{C}$ or $\mathbf{D}$, which contradicts part (ii) of the Theorem. \qed

\providecommand{\bysame}{\leavevmode\hbox to3em{\hrulefill}\thinspace}
\providecommand{\MR}{\relax\ifhmode\unskip\space\fi MR }
\providecommand{\MRhref}[2]{%
  \href{http://www.ams.org/mathscinet-getitem?mr=#1}{#2}
}
\providecommand{\href}[2]{#2}

\end{document}